# WIENER CHAOS SOLUTIONS OF LINEAR STOCHASTIC EVOLUTION EQUATIONS

By S. V. Lototsky[1] and B. L. Rozovskii[2]

*University of Southern California*

A new method is described for constructing a generalized solution of a stochastic evolution equation. Existence, uniqueness, regularity and a probabilistic representation of this Wiener Chaos solution are established for a large class of equations. As an application of the general theory, new results are obtained for several types of the passive scalar equation.

**1. Introduction.** Consider a stochastic evolution equation

$$(1.1) \qquad du(t) = (\mathcal{A}u(t) + f(t))\,dt + (\mathcal{M}u(t) + g(t))\,dW(t),$$

where $\mathcal{A}$ and $\mathcal{M}$ are differential operators, and $W$ is a cylindrical Brownian motion on a stochastic basis $\mathbb{F} = (\Omega, \mathcal{F}, \{\mathcal{F}_t\}_{t\geq 0}, \mathbb{P})$. Let $\mathcal{A}_0$ and $\mathcal{M}_0$ be the leading (highest-order) terms of the operators $\mathcal{A}$ and $\mathcal{M}$, respectively, and

$$\mathcal{Z}_t(u) = \int_0^t (\mathcal{M}u(s) + g(s))\,dW(s).$$

Traditionally, (1.1) is studied under the following assumptions:

 (i) The operator $\mathcal{A}_0 - \frac{1}{2}\mathcal{M}_0\mathcal{M}_0^*$ is elliptic.
 (ii) The noise term $\mathcal{Z}_t(v)$ is sufficiently regular. More specifically, for a suitable function space $\mathbf{X}$, $\mathbb{E}\int_0^T \|\mathcal{Z}_t(v)\|_{\mathbf{X}}^2\,dt < \infty$ for all $v$ in a dense subspace of $\mathbf{X}$.

Under these assumptions, there exists a unique Itô (strong) solution or a martingale (weak) solution $u$ of (1.1) so that $u \in L_2(\Omega \times (0,T); \mathbf{X})$ for $T > 0$

Received July 2004; revised December 2004.

[1]Supported by a Sloan Research Fellowship, NSF Career Award DMS-02-37724 and ARO Grant DAAD19-02-1-0374.

[2]Supported by ARO Grant DAAD19-02-1-0374 and ONR Grant N0014-03-1-0027.

(see, e.g., [5, 19, 33, 35]). In the future, we will refer to such solutions as *traditional* or *square integrable* solutions.

There are important examples demonstrating that the assumptions (i)–(ii) are necessary for the existence of a square integrable solution of (1.1).

In particular, it was shown in [30] that the stochastic advection–diffusion equation

$$\frac{\partial}{\partial t} u(t,x) = \Delta u(t,x) + u(t,x)\dot{W}(t,x), \qquad x \in \mathbb{R}^d,$$

where $\dot{W}(t,x)$ is the space–time white noise, has no square integrable solutions in any Sobolev space $\mathbf{X} = H^s(\mathbb{R}^d), s \in \mathbb{R}$, if $d \geq 2$. In this case assumption (i) holds but assumption (ii) does not.

Simple calculations show (see, e.g., [35]) that

$$du(t,x) = a^2 \frac{\partial^2}{\partial x^2} u(t,x)\,dt + \sigma \frac{\partial}{\partial x} u(t,x)\,dw(t),$$

where $w(t)$ is a one-dimensional Brownian motion and $u(0,x)$ is square integrable, has no square integrable solutions unless $a^2 - \frac{1}{2}\sigma^2 \geq 0$, that is, unless (i) holds. Obviously, (i) implies that the order of the operator $\mathcal{M}$ must be no larger than a half of the order of $\mathcal{A}$. Moreover, this assumption is, in some way, counterintuitive. Indeed, in the deterministic theory, only the highest-order operator ($\mathcal{A}_0$, in our case) has to be elliptic.

The objective of the current paper is to study stochastic differential equations of the type (1.1) without the restrictive assumptions (i) and (ii). The basic idea can be described as follows. If (1.1) does have a sufficiently regular solution, this solution can be projected on an orthonormal basis in some Hilbert space, resulting in a system of equations for the corresponding Fourier coefficients. We now turn this argument around, and *define* the solution of (1.1) as a formal Fourier series with the coefficients computed by solving the corresponding system. It often happens that this system has a solution under more general conditions than the original equation.

This approach could be traced to the classical separation of variables ideas in PDE's. For example, the Navier–Stokes equation is often defined as a system of coupled ODE's for the modes of its formal Fourier expansion with respect to the spatial variables (see, e.g., [7] and [26]). Similarly, in our case, the nonrandom spatiotemporal variables $(x,t)$ are being separated from the "random variable" (Brownian motion).

More specifically, the traditional solution of (1.1) is an $\mathcal{F}_T^W$-measurable square integrable random variable taking values in a Hilbert space $\mathbf{X}$. A classical result by Cameron and Martin [4] provides an orthonormal basis $\Xi = \{\xi_\alpha, \alpha \in \mathcal{J}\}$ in the Wiener Chaos space $L_2(\mathbb{W}; \mathbf{X})$, where $\mathcal{J}$ is the set of multi-indices $\alpha = \{\alpha_i^k\}$ of finite length $|\alpha| = \sum_{i,k=1}^\infty \alpha_i^k$. Accordingly, a *Wiener Chaos solution* of (1.1) is defined as a formal Fourier series with respect to the Cameron–Martin basis $\Xi$. By construction, this solution is strong in the probabilistic sense, that is, uniquely determined by the coefficients, free terms, initial condition and the Wiener process. The coefficients in the Fourier series are computed by solving the corresponding *propagator*, a lower-triangular system of deterministic parabolic equations, uniquely determined by (1.1) (in earlier works, e.g., [28], the propagator is referred to as S-system).

Of course, unless both (i) and (ii) hold, one could not expect the resulting Fourier series to converge in $L_2(\mathbb{W}; \mathbf{X})$. However, we demonstrate that under quite general assumptions, the Fourier series converges in a "minimal" weighted Wiener Chaos space $L_{2,Q}(\mathbb{W}; \mathbf{X})$. The construction of this space is quite simple (cf. [30]). Given a sequence of positive numbers $Q = \{q_1, q_2, \ldots\}$, we define $L_{2,Q}(\mathbb{W}; \mathbf{X})$ as the collection of sequences $u = \{u_\alpha, \alpha \in \mathcal{J}\}$ with $u_\alpha \in \mathbf{X}$ so that

$$\|u\|^2_{L_{2,Q}(\mathbb{W};\mathbf{X})} := \sum_{\alpha \in \mathcal{J}} q^{2\alpha} \|u_\alpha\|^2_{\mathbf{X}} < \infty,$$

where

$$q^\alpha = \prod_{i,k} q_k^{\alpha_i^k}, \qquad \alpha \in \mathcal{J}.$$

We remark that the space $L_2(\mathbb{W}; \mathbf{X})$ is a very special example of the sequence spaces studied in [15, 16, 17].

Of course, the Wiener Chaos solution is a weak solution designed for dealing with equations that do not have square integrable solutions. However, in some problems, the Wiener Chaos solution serves as a convenient first step in the investigation of square integrable solutions. For example, this is the case in the proof of the existence of square integrable solutions of degenerate parabolic SPDE's (Corollary 4.2 and Theorem 6.3).

Constructions based on various forms of the Wiener Chaos decomposition are popular in the study of stochastic differential equations, both ordinary and with partial derivatives. For stochastic ordinary differential equations, [20] used multiple Wiener integral expansion to study Itô's diffusions with nonsmooth coefficients. More recently, LeJan and Raimond [23] used a similar approach in the construction of stochastic flows. Various versions of the Wiener Chaos appear in a number of papers on nonlinear filtering and related topics; see, for example, [3, 24, 28, 31, 37]. The book by Holden et al. [12] presents a systematic approach to the stochastic differential equations based on the white noise theory; see also [11, 34] and the references therein.

The propagator was first introduced by Mikulevicius and Rozovskii [27], and further studied in [24], as a numerical tool for solving the Zakai filtering equation. In [30], the propagator was used to construct a generalized solution of the reaction–diffusion equation driven by the space–time white noise in several space dimensions. A similar system can be derived for certain nonlinear equations, such as the stochastic Navier–Stokes equation [29].

The main results of the paper are:

(1) Existence, uniqueness, regularity and the Krylov–Veretennikov formula for the Wiener Chaos solution of (1.1) (Theorems 3.4, 4.1 and Corollary 4.2).

(2) A Feynman–Kac formula for Wiener Chaos solutions in $L_{2,Q}(\mathbb{W}; L_2(\mathbb{R}^d))$ (Theorem 5.1).

(3) Existence, uniqueness and regularity properties for the transport equation with:

   (a) a space–time white noise-type velocity field (Theorem 6.1);
   (b) an incompressible Kraichnan turbulent velocity field, the nonviscous case (Theorem 6.3).

The paper is organized as follows. Section 2 introduces the Cameron–Martin basis and the weighted Wiener Chaos spaces. Section 3 presents the general definition of the Wiener Chaos solution and establishes the connection with the traditional and the white noise solutions. In the three main Sections, 4, 5 and 6, three types of stochastic equations are studied under three different sets of assumptions; these assumptions are listed at the beginning of the corresponding section. Section 4 presents the basic existence/uniqueness/regularity result for an abstract stochastic evolution equation under assumptions (A1)–(A3). Section 5 establishes probabilistic and multiple Wiener integral representations of Wiener Chaos solutions for stochastic partial differential equations in $\mathbb{R}^d$ with nonrandom coefficients under assumptions (B0)–(B5). Finally, Section 6 illustrates the results of the previous two sections for several versions of the turbulent transport equation under assumptions (S1)–(S3).

The following notation will be in force throughout the paper. $\Delta$ is the Laplace operator, $D_i = \partial/\partial x_i$, $i = 1, \ldots, d$, and summation over the repeated indices is assumed. The space of continuous functions is denoted by **C**, and $H_2^\gamma(\mathbb{R}^d)$, $\gamma \in \mathbb{R}$, is the Sobolev space

$$\left\{ f : \int_{\mathbb{R}^d} |\hat{f}(y)|^2 (1+|y|^2)^\gamma \, dy < \infty \right\}, \qquad \text{where } \hat{f} \text{ is the Fourier transform of } f.$$

**2. Weighted Wiener Chaos spaces.** In this section we review the construction of the Cameron–Martin basis and define the spaces of generalized random elements.

For a fixed $T > 0$, let $\mathbb{F} = (\Omega, \mathcal{F}, \{\mathcal{F}_t\}_{0 \leq t \leq T}, \mathbb{P})$ be a stochastic basis with the usual assumptions and let $W = (w_k = w_k(t), k \geq 1, 0 \leq t \leq T)$ be a collection of independent standard Brownian motions on $\mathbb{F}$. Denote by $\mathcal{F}_t^W$ the $\sigma$-algebra generated by the random variables $\{w_k(s), k \geq 1, s \leq t\}$, and by $L_2(\mathbb{W})$, the Hilbert space of $\mathcal{F}_T^W$-measurable square integrable random variables.

Our first step is to construct the Cameron–Martin basis, a special orthonormal basis in the space $L_2(\mathbb{W})$.

Let $\mathfrak{m} = \{m_k, k \geq 1\}$ be an orthonormal basis in $L_2((0,T))$ so that each function $m_k = m_k(t)$ is bounded for $t \in [0,T]$. Given such a basis $\mathfrak{m}$, define independent standard Gaussian random variables

$$\xi_{ik} = \int_0^T m_i(s)\, dw_k(s).$$

Consider the set of multi-indices

$$\mathcal{J} = \left\{ \alpha = (\alpha_i^k, i, k \geq 1), \alpha_i^k \in \{0, 1, 2, \ldots\}, \sum_{i,k} \alpha_i^k < \infty \right\}.$$

The set $\mathcal{J}$ is countable, and, for every $\alpha \in \mathcal{J}$, only finitely many of $\alpha_i^k$ are not equal to zero. For $\alpha \in \mathcal{J}$, we write

$$|\alpha| = \sum_{i,k} \alpha_i^k, \qquad \alpha! = \prod_{i,k} \alpha_i^k!,$$

and define the collection $\Xi = \{\xi_\alpha, \alpha \in \mathcal{J}\}$ of random variables so that

$$(2.1) \qquad \xi_\alpha = \frac{1}{\sqrt{\alpha!}} \prod_{i,k} H_{\alpha_i^k}(\xi_{ik}),$$

where

$$H_n(t) = (-1)^n e^{t^2/2} \frac{d^n}{dt^n} e^{-t^2/2}$$

is $n$th Hermite polynomial.

THEOREM 2.1. *The collection $\Xi = \{\xi_\alpha, \alpha \in \mathcal{J}\}$ is an orthonormal basis in $L_2(\mathbb{W})$.*

PROOF. This is a version of the classical result by Cameron and Martin [4]. □

Using the Cameron–Martin basis $\Xi$ in $L_2(\mathbb{W})$, we now define the space of generalized random elements.

Let $\mathbf{X}$ be a Banach space and $Q = \{q_1, q_2, \ldots\}$, a sequence of positive numbers. Define

$$(2.2) \qquad q^\alpha = \prod_{i,k} q_k^{\alpha_i^k}, \qquad \alpha \in \mathcal{J}.$$

DEFINITION 2.2. (i) The $Q$-weighted Wiener Chaos space $L_{2,Q}(\mathbb{W}; \mathbf{X})$ is the collection of sequences $u = \{u_\alpha, \alpha \in \mathcal{J}\}$ with $u_\alpha \in \mathbf{X}$ so that

$$\|u\|^2_{L_{2,Q}(\mathbb{W};\mathbf{X})} := \sum_{\alpha \in \mathcal{J}} q^{2\alpha} \|u_\alpha\|^2_{\mathbf{X}} < \infty.$$

(ii) The $Q^-$-weighted Wiener Chaos space $L_{2,Q^-}(\mathbb{W}; \mathbf{X})$ is the collection of sequences $u = \{u_\alpha, \alpha \in \mathcal{J}\}$ with $u_\alpha \in \mathbf{X}$ so that

$$\|u\|^2_{L_{2,Q^-}(\mathbb{W};\mathbf{X})} := \sum_{\alpha \in \mathcal{J}} q^{-2\alpha} \|u_\alpha\|^2_{\mathbf{X}} < \infty.$$

The space $\mathbf{X}$ will be omitted from the notation if $\mathbf{X} = \mathbb{R}$. The sequence $Q$ will be omitted from the notation if $Q = 1$, that is, if $q_k = 1$ for all $k$.

Given a $u = \{u_\alpha, \alpha \in \mathcal{J}\} \in L_{2,Q}(\mathbb{W}; \mathbf{X})$, we call each $u_\alpha$ a *generalized Fourier coefficient* of $u$ and identify $u$ with a formal Fourier series

$$u = \sum_{\alpha \in \mathcal{J}} u_\alpha \xi_\alpha.$$

The members of the set $L_{2,Q}(\mathbb{W}; \mathbf{X})$ are called $\mathbf{X}$-*valued generalized random elements*. Similarly, the members of the set $L_{2,Q}(\mathbb{W}; L_1((0,T); \mathbf{X}))$ are called $\mathbf{X}$-*valued generalized random processes*.

For $u \in L_{2,Q}(\mathbb{W}; \mathbf{X})$ and $v \in L_{2,Q^-}(\mathbb{W})$ we define

(2.3) $$\langle\!\langle u, v \rangle\!\rangle = \sum_{\alpha \in \mathcal{J}} u_\alpha v_\alpha;$$

the series in (2.3) converges in the norm of $\mathbf{X}$ by the Cauchy–Schwartz inequality.

**3. Linear stochastic evolution equations and the propagator.** In this section we define the Wiener Chaos solution and establish its connection with the traditional and white noise solutions. To motivate the definition, we start by reviewing the main results about the traditional solution.

Let $(V, H, V')$ be a normal triple of Hilbert spaces so that $V \subset H \subset V'$ with both embeddings continuous; for the complete definition of the normal triple see Section 3.1 in [35]. Denote by $\langle v', v \rangle$, $v' \in V'$, $v \in V$, the duality between $V$ and $V'$ relative to the inner product in $H$.

For $t \in [0,T]$, consider families of linear operators $\mathcal{A} = \mathcal{A}(t)$ and $\mathcal{M}_k = \mathcal{M}_k(t)$ so that, for each $t$, the operators $\mathcal{A}(t): V \to V'$, $\mathcal{M}_k(t): V \to H$ are bounded. Consider the following equation:

(3.1) $$u(t) = u_0 + \int_0^t (\mathcal{A}u(s) + f(s))\,ds + \int_0^t (\mathcal{M}_k u(s) + g_k(s))\,dw_k(s),$$

$$0 \le t \le T.$$

Recall that summation convention over the pairs of repeated indices is in force.

We proceed with a review of the traditional approach. Assume that, for all $v \in V$, $t \in [0, T]$,

$$\sum_{k \geq 1} \|\mathcal{M}_k(t) v\|_H^2 < \infty, \tag{3.2}$$

and the nonrandom input data $u_0, f$ and $g_k$ satisfy

$$\|u_0\|_H^2 + \int_0^T \|f(t)\|_{V'}^2 \, dt + \sum_{k \geq 1} \int_0^T \|g_k(t)\|_H^2 \, dt < \infty. \tag{3.3}$$

DEFINITION 3.1. An $\mathcal{F}_t^W$-adapted process $u \in L_2(\mathbb{W}; L_2((0,T); V))$ is called a *square integrable solution* of (3.1) if, for every $v \in V$, there exists a measurable subset $\Omega'$ of $\Omega$ with $\mathbb{P}(\Omega') = 1$, so that, for all $0 \leq t \leq T$, the equality

$$\begin{aligned}(u(t), v)_H &= (u_0, v)_H + \int_0^t \langle \mathcal{A} u(s) + f(s), v \rangle \, ds \\ &\quad + \sum_{k \geq 1} \int_0^t (\mathcal{M}_k u(s) + g_k(s), v)_H \, dw_k(s)\end{aligned} \tag{3.4}$$

holds on $\Omega'$. Similarly, $u \in L_2((0,T); V)$ is a *solution* of the deterministic equation

$$u(t) = u_0 + \int_0^t \mathcal{A} u(s) \, ds + \int_0^t f(s) \, ds$$

if, for every $v \in V$ and $t \in [0, T]$, equality (3.4) holds with $\mathcal{M}_k = g_k = 0$.

Existence and uniqueness of the traditional solution of (3.1) are established under an additional assumption about the operators $\mathcal{A}$ and $\mathcal{M}_k$.

DEFINITION 3.2. Equation (3.1) is called *strongly parabolic* if there exists a positive number $\varepsilon$ and a real number $C_0$ so that, for all $v \in V$ and $t \in [0, T]$,

$$2 \langle \mathcal{A}(t) v, v \rangle + \sum_{k \geq 1} \|\mathcal{M}_k(t) v\|_H^2 + \varepsilon \|v\|_V^2 \leq C_0 \|v\|_H^2. \tag{3.5}$$

Equation (3.1) is called *weakly parabolic* (or degenerate parabolic) if condition (3.5) holds with $\varepsilon = 0$.

THEOREM 3.3. *If* (3.3) *holds and* (3.1) *is strongly parabolic, then there exists a unique square integrable solution of* (3.1). *The solution process* u *belongs to*

$$L_2(\mathbb{W}; L_2((0,T); V)) \cap L_2(\mathbb{W}; \mathbf{C}((0,T), H))$$

*and satisfies*

$$\mathbb{E}\bigg(\sup_{0\leq t\leq T} \|u(t)\|_H^2 + \int_0^T \|u(t)\|_V^2 \, dt\bigg)$$
(3.6)
$$\leq C(C_0, \varepsilon, T)\bigg(\|u_0\|_H^2 + \int_0^T \|f(t)\|_{V'}^2 \, dt + \sum_{k\geq 1} \int_0^T \|g_k(t)\|_H^2 \, dt\bigg).$$

PROOF. This follows, for example, from Theorem 3.1.4 in [35]. □

When (3.1) is weakly parabolic, then the solvability result is somewhat different; see Section 3.2 in [35] for details.

As an element of the Hilbert space $L_2(\mathbb{W}; L_2((0,T); V))$, the traditional solution of (3.1) admits a representation $u(t) = \sum_{\alpha\in\mathcal{J}} u_\alpha(t)\xi_\alpha$ in the Cameron–Martin basis $\Xi$.

THEOREM 3.4. *An* $\mathcal{F}_t^W$*-adapted process* u *is a square integrable solution of* (3.1) *if and only if* $u(t) = \sum_{\alpha\in\mathcal{J}} u_\alpha(t)\xi_\alpha$ *so that the Fourier coefficients* $u_\alpha$ *satisfy*

(3.7) $$\sum_{\alpha\in\mathcal{J}} \bigg(\int_0^T \|u_\alpha(t)\|_V^2 \, dt + \sup_{0\leq t\leq T} \|u_\alpha(t)\|_H^2\bigg) < \infty,$$

*and solve the* propagator

$$u_\alpha(t) = u_0 I(|\alpha| = 0) + \int_0^t (\mathcal{A}u_\alpha(s) + f(s)I(|\alpha| = 0)) \, ds$$
(3.8)
$$+ \int_0^t \sum_{i,k} \sqrt{\alpha_i^k} (\mathcal{M}_k u_{\alpha^-(i,k)}(s) + g_k(s)I(|\alpha| = 1)) m_i(s) \, ds,$$

*where* $\alpha^-(i,k)$ *is the multi-index with components*

$$(\alpha^-(i,k))_j^l = \begin{cases} \max(\alpha_i^k - 1, 0), & \text{if } i = j \text{ and } k = l, \\ \alpha_j^l, & \text{otherwise.} \end{cases}$$

Before presenting the proof, we define the Wiener Chaos solution of (3.1). The definition is motivated by Theorem 3.4.

DEFINITION 3.5. A $V$-valued generalized random process $u$ is called a *Wiener Chaos solution* of (3.1) if the generalized Fourier coefficients $u_\alpha, \alpha \in \mathcal{J}$, of $u$ are a solution of the propagator (3.8).

To prove Theorem 3.4 and to derive an alternative characterization of the Wiener Chaos solution, we need a few additional constructions.

Denote by $\mathcal{H}$ the set

$$\mathcal{H} = \bigcup_{n \geq 1} L_\infty((0,T); \mathbb{R}^n).$$

If $h \in \mathcal{H}$, then there exists an $N \geq 1$ so that

$$h = (h_1, \ldots, h_N),$$

with each $h_k \in L_\infty((0,T))$. We define

$$\mathcal{E}(t,h) = \exp\left(\sum_{k=1}^N \left(\int_0^t h_k(s)\, dw_k(s) - \tfrac{1}{2}\int_0^t |h_k(s)|^2\, ds\right)\right), \qquad h \in \mathcal{H},$$

(3.9) $\quad \mathcal{E}(h) = \mathcal{E}(T,h),$

$$h_{ik} = \int_0^T h_k(t) m_i(t)\, dt, \qquad m_i \in \mathfrak{m},$$

and

$$h^\alpha = \prod_{i,k} h_{ik}^{\alpha_i^k}, \qquad \alpha \in \mathcal{J}.$$

The following properties of the process $\mathcal{E}(t,h)$ are verified by direct calculation:

(3.10) $$\mathcal{E}(h) = \sum_{\alpha \in \mathcal{J}} \frac{h^\alpha}{\sqrt{\alpha!}},$$

(3.11) $$\mathcal{E}(t,h) = 1 + \int_0^t \mathcal{E}(s,h) h_k(s)\, dw_k(s).$$

One consequence of (3.10) is that $\mathcal{E}(h) \in L_{2,Q}(\mathbb{W})$ for every weight sequence $Q$, with $\|\mathcal{E}(h)\|_{L_{2,Q}(\mathbb{W})}^2 = \exp(\sum_{k=1}^N q_k^2 \|h_k\|_{L_2((0,T))}^2)$. Another consequence of (3.10) is an alternative representation of $\xi_\alpha$:

(3.12) $$\xi_\alpha = \frac{1}{\sqrt{\alpha!}} \frac{\partial^{|\alpha|}}{\partial h^\alpha} \mathcal{E}(h)\bigg|_{h=0}.$$

REMARK 3.6. By Lemma 4.3.2 in [32], the family $\{\mathcal{E}(h), h \in \mathcal{H}\}$ is dense in $L_2(\mathbb{W})$ and therefore in every $L_{2,Q}(\mathbb{W})$.

PROOF OF THEOREM 3.4. (i) Assume that $u = u(t)$ is a traditional solution of (3.1). Equation (3.12) implies

$$u_\alpha(t) = \mathbb{E}(u(t)\xi_\alpha) = \frac{1}{\sqrt{\alpha!}} \frac{\partial^{|\alpha|}}{\partial h^\alpha} \mathbb{E}(u(t)\mathcal{E}(h))\bigg|_{h=0}.$$

Using the $\mathcal{F}_t^W$-measurability of $u(t)$ and the martingale property (3.11) of $\mathcal{E}(t,h)$, we derive

$$\mathbb{E}(u(t)\mathcal{E}(h)) = \mathbb{E}(u(t)\mathbb{E}(\mathcal{E}(h)|\mathcal{F}_t^W)) = \mathbb{E}(u(t)\mathcal{E}(t,h)),$$

and (3.8) follows after applying the Itô formula to the product $u(t)\mathcal{E}(t,h)$ and differentiating the resulting equation with respect to $h$. By Theorem 3.3.3(iii) in [22], this differentiation with respect to $h$ is justified.

The same arguments show that the time evolution of

(3.13) $$\xi_\alpha(t) := \mathbb{E}(\xi_\alpha|\mathcal{F}_t^W)$$

is described by

(3.14) $$\xi_\alpha(t) = I(|\alpha| = 0) + \int_0^t \sum_{i,k} \sqrt{\alpha_i^k} \xi_{\alpha^-(i,k)}(s) m_i(s) \, dw_k(s).$$

(ii) Conversely, assume that condition (3.7) holds. Then the process $u(t) = \sum_{\alpha \in \mathcal{J}} u_\alpha(t)\xi_\alpha$ satisfies

$$u \in L_2(\mathbb{W}; L_2((0,T); V)) \cap L_2(\mathbb{W}; \mathbf{C}((0,T); H)),$$

and the function $u_h = \mathbb{E}(u\mathcal{E}(h))$ belongs to $L_2((0,T); V) \cap \mathbf{C}((0,T); H)$. By Parseval's identity and relation (3.10),

$$u_h = \sum_{\alpha \in \mathcal{J}} \frac{h^\alpha}{\sqrt{\alpha!}} u_\alpha.$$

Using (3.8), we conclude that, for every $\varphi \in V$ and $h \in \mathcal{H}$, the function $u_h$ satisfies

$$(u_h(t), \varphi)_H = (u_0, \varphi)_H + \int_0^t (\mathcal{A}u_h(s), \varphi)_H \, ds + \int_0^t \langle f(s), \varphi \rangle \, ds$$
$$+ \sum_{\alpha \in \mathcal{J}} \frac{h^\alpha}{\alpha!} \sum_{i,k} \int_0^t \sqrt{\alpha_i^k} m_i(s)((\mathcal{M}_k u_{\alpha^-(i,k)}(s), \varphi)_H$$
$$+ (g_k(s), \varphi)_H I(|\alpha| = 1)) \, ds.$$

To simplify the above equality note that if $I(t) = \int_0^t (\mathcal{M}_k u(s), \varphi)_H \, dw_k(s)$, then the $\mathcal{F}_t^W$-measurability of $I(t)$ and relation (3.14) imply

(3.15) $$\mathbb{E}(I(t)\xi_\alpha) = \mathbb{E}(I(t)\xi_\alpha(t)) = \int_0^t \sum_{i,k} \sqrt{\alpha_i^k} m_i(s)(\mathcal{M}_k u_{\alpha^-(i,k)}(s), \varphi)_H \, ds.$$

WIENER CHAOS SOLUTIONS  11

Similarly,

$$\mathbb{E}\left(\xi_\alpha \int_0^t g_k(s)\, dw_k(s)\right) = \sum_{i,k} \int_0^t \sqrt{\alpha_i^k}\, m_i(s)(g_k(s),\varphi)_H I(|\alpha|=1)\, ds.$$

Therefore,

$$\sum_{\alpha\in\mathcal{J}} \frac{h^\alpha}{\alpha!} \sum_{i,k} \int_0^t \sqrt{\alpha_i^k}\, m_i(s)((\mathcal{M}_k u_{\alpha^-(i,k)}(s),\varphi)_H + (g_k(s),\varphi)_H I(|\alpha|=1))\, ds$$
$$= \mathbb{E}\left(\mathcal{E}(h) \int_0^t ((\mathcal{M}_k u(s),\varphi)_H + (g_k(s),\varphi)_H)\, dw_k(s)\right).$$

As a result,

$$\mathbb{E}(\mathcal{E}(h)(u(t),\varphi)_H) = \mathbb{E}(\mathcal{E}(h)(u_0,\varphi)_H) + \mathbb{E}\left(\mathcal{E}(h)\int_0^t \langle \mathcal{A}u(s),\varphi\rangle\, ds\right)$$

(3.16)
$$+ \mathbb{E}\left(\mathcal{E}(h)\int_0^t \langle f(s),\varphi\rangle\, ds\right)$$
$$+ \mathbb{E}\left(\mathcal{E}(h)\int_0^t ((\mathcal{M}_k u(s),\varphi)_H + (g_k(s),\varphi)_H)\, dw_k(s)\right).$$

Equality (3.16) and Remark 3.6 imply that, for each $t$ and each $\varphi$, (3.4) holds with probability 1. Due to continuity of $u$, a single probability-1 set can be chosen for all $t \in [0,T]$.

Theorem 3.4 is proved. □

Theorem 3.4 implies that, if it exists, the traditional solution of (3.1) coincides with the Wiener Chaos solution.

We now give an alternative characterization of the Wiener Chaos solution.

THEOREM 3.7. *A $V$-valued generalized random process $u$ is a Wiener Chaos solution of (3.1) if and only if, for every $h \in \mathcal{H}$, the function $u_h(t) = \langle\!\langle u(t),\mathcal{E}(h)\rangle\!\rangle$ is a solution of the equation*

$$(3.17)\quad u_h(t) = u_0 + \int_0^t (\mathcal{A}u_h(s) + f(s) + h_k(s)\mathcal{M}_k u_h(s) + h_k(s)g_k(s))\, ds.$$

PROOF. Assume that $u$ is a Wiener Chaos solution of (3.1), that is, $u \in L_{2,Q}(\mathbb{W};V)$ for some $Q$ and each $u_\alpha$ is a solution of (3.8). By definition (2.3) and relation (3.10),

$$(3.18)\qquad\qquad\qquad u_h = \sum_{\alpha\in\mathcal{J}} \frac{h^\alpha}{\sqrt{\alpha!}} u_\alpha.$$

Then (3.17) follows from (3.8). Indeed,

$$\sum_\alpha \frac{h^\alpha}{\sqrt{\alpha!}} \sum_{i,k} \sqrt{\alpha_i^k} \mathcal{M}_k u_{\alpha^-(i,k)} m_i = \sum_\alpha \sum_{i,k} \frac{h^{\alpha^-(i,k)}}{\sqrt{\alpha^-(i,k)!}} \mathcal{M}_k u_{\alpha^-(i,k)} m_i h_{ik}$$
$$= \sum_{i,k} \left( \sum_\alpha \frac{h^\alpha}{\sqrt{\alpha!}} \mathcal{M}_k u_\alpha \right) m_i h_{ik} = h_k \mathcal{M}_k u_h.$$

Computations for the other terms are similar.

Conversely, assume that $u \in L_{2,Q}(\mathbb{W}; L_2((0,T); V))$ and $u_h \in L_2((0,T); V)$ is a solution of (3.17). By relation (3.18),

$$u_\alpha = \frac{1}{\sqrt{\alpha!}} \frac{\partial^{|\alpha|}}{\partial h^\alpha} u_h \bigg|_{h=0}.$$

Then term-by-term differentiation in (3.17) implies (3.8). Theorem 3.7 is proved. $\square$

To conclude this section, we briefly discuss the relation between the Wiener Chaos and white noise solutions. While the white noise solution for (3.1) has not been defined in general, analysis of the particular cases [11, 12, 34] shows that the white noise solution exists provided (3.17) has a solution that is an *analytic* function of $h$. The white noise solution of (3.1) is constructed as an element of the space

$$(\mathcal{S})_{-\rho}(\mathbf{X}) = \left\{ u_\alpha : \sum_{\alpha \in \mathcal{J}} r_{\alpha,\ell}^2 \|u_\alpha\|_\mathbf{X}^2 < \infty \text{ for some } \ell > 0 \right\},$$

where

$$r_{\alpha,\ell}^2 = (\alpha!)^{-\rho} \prod_{i,k} (2ik)^{-\ell \alpha_i^k}, \qquad \rho \in [0,1];$$

$(\mathcal{S})_{-0}(\mathbb{R}^d)$ is known as the space of Hida distributions [11]. Even though $(\mathcal{S})_{-\rho}(\mathbf{X})$ is not, in general, related to $L_{2,Q}(\mathbb{W}; \mathbf{X})$, in many examples one can take either $q_k = q$ or $q_k = Cq^k$ for some $q < 1$. For such weights,

$$L_{2,Q}(\mathbb{W}; \mathbf{X}) \subset (\mathcal{S})_{-0}(\mathbf{X})$$

with strict inclusion, so that the Wiener Chaos approach provides better regularity results. In addition, we remark that in contrast to the Wiener Chaos solutions, the white noise solutions are weak solutions in the probabilistic sense.

**4. Linear evolution equations in weighted Wiener Chaos spaces.** In this section we present the main existence, uniqueness and regularity result for the Wiener Chaos solution of (3.1). We make the following assumptions:

(A1) There exist positive numbers $C_1$ and $\delta$ so that, for every $v \in V$ and $t \in [0,T]$,

(4.1) $$\langle \mathcal{A}(t)v, v\rangle + \delta\|v\|_V^2 \leq C_1\|v\|_H^2.$$

(A2) There exists a real number $C_2$ and a sequence of positive numbers $Q = \{q_k, k \geq 1\}$ so that, for every $v \in H$ and $t \in [0,T]$,

(4.2) $$2\langle \mathcal{A}(t)v, v\rangle + \sum_{k\geq 1} q_k^2 \|\mathcal{M}_k(t)v\|_H^2 \leq C_2\|v\|_H^2.$$

(A3) The initial condition $u_0$ is nonrandom and belongs to $H$; the process $f = f(t)$ is deterministic and $\int_0^T \|f(t)\|_{V'}^2 \, dt < \infty$; each $g_k = g_k(t)$ is a deterministic process and

$$\sum_{k\geq 1} \int_0^T q_k^2 \|g_k(t)\|_H^2 \, dt < \infty.$$

Denote by $(P_{t,s}, 0 \leq s \leq t \leq T)$ the semigroup generated by the operator $\mathcal{A}$; $P_t := P_{t,0}$. In other words, the solution of the equation

$$v(t) = v_0 + \int_s^t \mathcal{A}v(\tau)\,d\tau + \int_s^t f(\tau)\,d\tau, \qquad 0 \leq s \leq t \leq T,$$

$v_0 \in H, f \in L_2((0,T); V')$, is written as

$$v(t) = P_{t,s}v_0 + \int_s^t P_{t,\tau}f(\tau)\,d\tau.$$

This solution exists by assumption (A1) and Theorem 3.3. Denote by $u_{(0)}$ the solution of (3.8) when $|\alpha| = 0$.

THEOREM 4.1. *Under assumptions* (A1)–(A3), (3.1) *has a unique Wiener Chaos solution. The solution $u = u(t)$ has the following properties:*

(1) *For every $\gamma \in (0,1)$,*

(4.3) $$u \in L_{2,\gamma Q}(\mathbb{W}; L_2((0,T); V)) \cap L_{2,\gamma Q}(\mathbb{W}; \mathbf{C}((0,T); H)).$$

(2) *For every $0 \leq t \leq T$, $u(t) \in L_{2,Q}(\mathbb{W}; H)$ and the following relations hold:*

$$\sum_{|\alpha|=n} q^\alpha u_\alpha(t)\xi_\alpha$$

$$= \sum_{k_1,\ldots,k_n \geq 1} \int_0^t \int_0^{s_n}$$

$$(4.4) \quad \cdots \int_0^{s_2} P_{t,s_n} \overline{\mathcal{M}}_{k_n} \cdots P_{s_2,s_1}$$
$$\times (\overline{\mathcal{M}}_{k_1} u_{(0)} + q_{k_1} g_{k_1}(s_1))\, dw_{k_1}(s_1) \cdots dw_{k_n}(s_n),$$
$$n \geq 1,$$

where $\overline{\mathcal{M}}_k = q_k \mathcal{M}_k$;

$$\|u(t)\|^2_{L_{2,Q}(\mathbb{W};H)}$$
$$(4.5) \quad \leq 3 e^{C_2 t} \left( \|u_0\|_H^2 + C_f \int_0^t \|f(s)\|_{V'}^2\, ds + \sum_{k\geq 1} \int_0^t q_k^2 \|g_k(s)\|_H^2\, ds \right),$$

where the number $C_2$ is from (4.2) and the positive number $C_f$ depends only on $\delta$ and $C_1$ from (4.1).

PROOF. The arguments build on the techniques from [25].

(1) By Theorem 3.3, there exists a unique traditional solution of equation

$$(4.6) \quad \begin{aligned} v(t) &= u_0 + \int_0^t (\mathcal{A}v(s) + f(s))\, ds \\ &\quad + \sum_{k\geq 1} \int_0^t (\mathcal{M}_k v(s) + g_k(s)) \gamma q_k\, dw_k(s). \end{aligned}$$

By Theorems 3.4 and 3.7, $v_\alpha = \gamma^{|\alpha|} q^\alpha u_\alpha$, and (4.3) follows.

(2) Direct computations show that, for $|\alpha| = 1$ with $\alpha_i^k = 1$, the corresponding solution $u_\alpha := u_{(ik)}$ of (3.8) is

$$(4.7) \quad u_{(ik)}(t) = \int_0^t P_{t,s}(\mathcal{M}_k u_{(0)}(s) + g_k(s)) m_i(s)\, ds.$$

Since $\xi_{(ik)} = \xi_{ik} = \int_0^T m_i(s)\, dw_k(s)$, we conclude that

$$\sum_{i\geq 1} u_{ik}(t) \xi_{(ik)} = \int_0^t P_{t,s}(\mathcal{M}_k u_{(0)}(s) + g_k(s))\, dw_k(s)$$

or

$$\sum_{|\alpha|=1} q^\alpha u_\alpha(t) \xi_\alpha = \sum_{k\geq 1} \int_0^t P_{t,s}(\overline{\mathcal{M}}_k u_{(0)}(s) + q_k g_k(s))\, dw_k(s).$$

Continuing by induction on $|\alpha|$ and using the relation between the Hermite polynomials and the iterated Itô integrals ([13], Theorem 3.1), we derive (4.4).

An immediate consequence of (4.4) is the following energy equality:

$$\sum_{|\alpha|=n} q^{2\alpha} \|u_\alpha(t)\|_H^2$$

$$(4.8) \quad = \sum_{k_1,\ldots,k_n \geq 1} \int_0^t \int_0^{s_n} \cdots \int_0^{s_2} \|P_{t,s_n}\overline{\mathcal{M}}_{k_n} \cdots P_{s_2,s_1}(\overline{\mathcal{M}}_{k_1} u_{(0)} + q_{k_1} g_{k_1})\|_H^2 \, ds^n,$$

where $ds^n = ds_1 \cdots ds_n$.

To derive (4.5), consider first the homogeneous equation with $f = g_k = 0$ and define $F_n(t) = \sum_{|\alpha|=n} q^{2\alpha} \|u_\alpha(t)\|_H^2$, $n \geq 0$. Note that $u_{(0)}(t) = P_t u_0$.

By assumption (4.2),

$$(4.9) \quad \frac{d}{dt} F_0(t) \leq C_2 F_0(t) - \sum_{k \geq 1} \|\overline{\mathcal{M}}_k P_t u_0\|_H^2.$$

For $n \geq 1$, the energy equality (4.8) implies

$$\frac{d}{dt} F_n(t) = \sum_{k_1,\ldots,k_n \geq 1} \int_0^t \int_0^{s_{n-1}} \cdots \int_0^{s_2} \|\overline{\mathcal{M}}_{k_n} P_{t,s_{n-1}} \overline{\mathcal{M}}_{k_{n-1}} \cdots \overline{\mathcal{M}}_{k_1} P_{s_1} u_0\|_H^2 \, ds^{n-1}$$

$$(4.10) \quad + \sum_{k_1,\ldots,k_n \geq 1} \int_0^t \int_0^{s_n} \cdots \int_0^{s_2} 2\langle \mathcal{A} P_{t,s_n} \overline{\mathcal{M}}_{k_n} \cdots \overline{\mathcal{M}}_{k_1} P_{s_1} u_0,$$

$$P_{t,s_n} \overline{\mathcal{M}}_{k_n} \cdots \overline{\mathcal{M}}_{k_1} P_{s_1} u_0 \rangle \, ds^n.$$

By assumption (4.2),

$$\sum_{k_1,\ldots,k_n \geq 1} \int_0^t \int_0^{s_n} \cdots \int_0^{s_2} 2\langle \mathcal{A} P_{t,s_n} \overline{\mathcal{M}}_{k_n} \cdots \overline{\mathcal{M}}_{k_1} P_{s_1} u_0,$$

$$(4.11) \quad P_{t,s_n} \overline{\mathcal{M}}_{k_n} \cdots \overline{\mathcal{M}}_{k_1} P_{s_1} u_0 \rangle \, ds^n$$

$$\leq - \sum_{k_1,\ldots,k_{n+1} \geq 1} \int_0^t \int_0^{s_n} \cdots \int_0^{s_2} \|\overline{\mathcal{M}}_{k_{n+1}} P_{t,s_n} \overline{\mathcal{M}}_{k_n} \cdots \overline{\mathcal{M}}_{k_1} P_{s_1} u_0\|_H^2 \, ds^n$$

$$+ C_2 \sum_{k_1,\ldots,k_n \geq 1} \int_0^t \int_0^{s_n} \cdots \int_0^{s_2} \|P_{t,s_n} \overline{\mathcal{M}}_{k_n} \cdots \overline{\mathcal{M}}_{k_1} P_{s_1} u_0\|_H^2 \, ds^n.$$

As a result, for $n \geq 1$,

$$\frac{d}{dt} F_n(t) \leq C_2 F_n(t)$$

$$(4.12) + \sum_{k_1,\ldots,k_n \geq 1} \int_0^t \int_0^{s_{n-1}} \cdots \int_0^{s_2} \|\overline{\mathcal{M}}_{k_n} P_{t,s_{n-1}} \overline{\mathcal{M}}_{k_{n-1}} \cdots \overline{\mathcal{M}}_{k_1} P_{s_1} u_0\|_H^2 \, ds^{n-1}$$

$$- \sum_{k_1,\ldots,k_{n+1}\geq 1} \int_0^t \int_0^{s_n} \cdots \int_0^{s_2} \|\overline{\mathcal{M}}_{k_{n+1}} P_{t,s_n} \overline{\mathcal{M}}_{k_n} \cdots \overline{\mathcal{M}}_{k_1} P_{s_1} u_0\|_H^2\, ds^n.$$

Consequently,

$$(4.13) \qquad \frac{d}{dt} \sum_{n=0}^{N} \sum_{|\alpha|=n} q^{2\alpha} \|u_\alpha(t)\|_H^2 \leq C_2 \sum_{n=0}^{N} \sum_{|\alpha|=n} q^{2\alpha} \|u_\alpha(t)\|_H^2,$$

so that, by the Gronwall inequality,

$$(4.14) \qquad \sum_{n=0}^{N} \sum_{|\alpha|=n} q^{2\alpha} \|u_\alpha(t)\|_H^2 \leq e^{C_2 t} \|u_0\|_H^2$$

or

$$(4.15) \qquad \|u(t)\|_{L_{2,Q}(\mathbb{W};H)}^2 \leq e^{C_2 t} \|u_0\|_H^2.$$

The remaining two cases, namely, $u_0 = g_k = 0$ and $u_0 = f = 0$, are analyzed in the same way, and then (4.5) follows by the triangle inequality. Theorem 4.1 is proved. $\square$

COROLLARY 4.2.  *Let $a_{ij}, b_i, c, \sigma_{ik}, \nu_k$ be deterministic measurable functions of $(t,x)$ so that*

$$|a_{ij}(t,x)| + |b_i(t,x)| + |c(t,x)| + |\sigma_{ik}(t,x)| + |\nu_k(t,x)| \leq K,$$

$i,j=1,\ldots,d, k\geq 1, x\in\mathbb{R}^d, 0\leq t\leq T;$

$$(a_{ij}(t,x) - \tfrac{1}{2}\sigma_{ik}(t,x)\sigma_{jk}(t,x))y_i y_j \geq 0,$$

$x,y\in\mathbb{R}^d, 0\leq t\leq T;$ *and*

$$\sum_{k\geq 1} |\nu_k(t,x)|^2 \leq C_\nu < \infty,$$

$x\in\mathbb{R}^d, 0\leq t\leq T.$ *Consider the equation*

$$(4.16) \qquad \begin{aligned} du &= (D_i(a_{ij} D_j u) + b_i D_i u + cu + f)\, dt \\ &\quad + (\sigma_{ik} D_i u + \nu_k u + g_k)\, dw_k. \end{aligned}$$

*Assume that the input data satisfy $u_0 \in L_2(\mathbb{R}^d)$, $f \in L_2((0,T); H_2^{-1}(\mathbb{R}^d))$, $\sum_{k\geq 1} \|g_k\|_{L_2((0,T)\times\mathbb{R}^d)}^2 < \infty$, and there exists an $\varepsilon > 0$ so that*

$$a_{ij}(t,x) y_i y_j \geq \varepsilon |y|^2, \qquad x,y \in \mathbb{R}^d, 0\leq t\leq T.$$

*Then there exists a unique Wiener Chaos solution $u = u(t,x)$ of (4.16). The solution has the following regularity:*

$$(4.17) \qquad u(t,\cdot) \in L_2(\mathbb{W}; L_2(\mathbb{R}^d)), \qquad 0 \leq t \leq T,$$

*and*

$$\mathbb{E}\|u\|^2_{L_2(\mathbb{R}^d)}(t)$$

(4.18)
$$\leq C^*\left(\|u_0\|^2_{L_2(\mathbb{R}^d)} + \|f\|^2_{L_2((0,T);H_2^{-1}(\mathbb{R}^d))} + \sum_{k\geq 1}\|g_k\|^2_{L_2((0,T)\times\mathbb{R}^d)}\right),$$

*where the positive number $C^*$ depends only on $C_\nu, K, T$ and $\varepsilon$.*

PROOF. Direct computation shows that the operators

$$\mathcal{A} = D_i(a_{ij}D_j) + b_iD_i + c, \qquad \mathcal{M}_k = \sigma_{ik}D_i + \nu_k$$

satisfy assumptions (A1) and (A2) in the normal triple $(H_2^1(\mathbb{R}^d), L_2(\mathbb{R}^d), H_2^{-1}(\mathbb{R}^d))$ with $q_k = 1$. Then relations (4.17) and (4.18) follow from (4.5). □

Note that in contrast to the statement of this corollary, all previous results concerning (4.16) required additional regularity of the coefficients and input data; see, for example, [35], Section 4.2.

Theorem 4.1 is a bona fide extension of Theorem 3.3. Indeed, if condition (3.5) holds so that (3.1) is strongly parabolic, then, taking $Q = 1$, we recover the statement of Theorem 3.3. Further analysis of condition (3.5) indicates that, for a strongly parabolic equation, one can find an admissible weight sequence $Q$ so that $q_k = q > 1$. Since for $Q > 1$ we have a strict inclusion $L_{2,Q}(\mathbb{W};\mathbf{X}) \subset L_2(\mathbb{W};\mathbf{X})$, Theorem 4.1 represents an improvement of Theorem 3.3 for strongly parabolic equations.

For weakly parabolic equations, similarly to the proof of Corollary 4.2, one can take $Q = 1$, and then Theorem 4.1 represents an extension of the existing results ([35], Sections 3.2 and 4.2).

For nonparabolic equations, the results of Theorem 4.1 are completely new.

REMARK 4.3. The formal Fourier series $u = \sum_{\alpha \in \mathcal{J}} u_\alpha \xi_\alpha$ for the Wiener Chaos solution is a generalization of the representation formula for solutions of SODE's, derived by ([20], Theorem 4) using iterated Itô integrals; see also [23], Theorem 3.2. Indeed,

$$\sum_{\alpha \in \mathcal{J}} u_\alpha \xi_\alpha = \sum_{n\geq 0}\left(\sum_{|\alpha|=n} u_\alpha \xi_\alpha\right),$$

and equality (4.4) connects the inner sum on the right with the iterated integrals. See also Example 5.2 below.

**5. Probabilistic representation of Wiener Chaos solutions.** As before, let $\mathbb{F} = (\Omega, \mathcal{F}, \{\mathcal{F}_t\}_{0 \le t \le T}, \mathbb{P})$ be a stochastic basis with the usual assumptions and let $W = \{w_k(t), k \ge 1, 0 \le t \le T\}$ be a collection of standard Wiener processes on $\mathbb{F}$.

Consider the linear equation in $\mathbb{R}^d$

$$(5.1) \quad du = (a_{ij} D_i D_j u + b_i D_i u + cu + f) \, dt + (\sigma_{ik} D_i u + \nu_k u + g_k) \, dw_k$$

under the following assumptions:

(B0) All coefficients, free terms and the initial condition are nonrandom.
(B1) The functions $a_{ij} = a_{ij}(t,x)$ are measurable and bounded in $(t,x)$, and
  (i)
  $$|a_{ij}(t,x) - a_{ij}(t,y)| \le C|x-y|, \qquad x,y \in \mathbb{R}^d, 0 \le t \le T,$$
  with $C$ independent of $t,x,y$;
  (ii) the matrix $(a_{ij})$ is uniformly positive definite, that is, there exists a $\delta > 0$ so that, for all vectors $y \in \mathbb{R}^d$ and all $(t,x)$, $a_{ij} y_i y_j \ge \delta |y|^2$.
(B2) The functions $b_i = b_i(t,x)$, $c = c(t,x)$ and $\nu_k = \nu_k(t,x)$ are measurable and bounded in $(t,x)$.
(B3) The functions $\sigma_{ik} = \sigma_{ik}(t,x)$ are continuous and bounded in $(t,x)$.
(B4) There exists a $p \ge d+1$ so that the functions $f = f(t,x)$ and $g_k = g_k(t,x)$ belong to
$$L_p((0,T); L_2(\mathbb{R}^d) \cap L_p(\mathbb{R}^d)).$$
(B5) The initial condition $u_0 = u_0(x)$ belongs to $L_2(\mathbb{R}^d) \cap W_p^2(\mathbb{R}^d)$, $p \ge d+1$, where $W_p^2$ is the Sobolev space $\{f : f, D_i f, D_i D_j f \in L_p(\mathbb{R}^d)\}$.

Under assumptions (B2)–(B4), there exists a sequence $Q = \{q_k, k \ge 1\}$ of positive numbers with the following properties:

(P1) The matrix $A$ with $A_{ij} = a_{ij} - (1/2) \sum_{k \ge 1} q_k \sigma_{ik} \sigma_{jk}$ satisfies $A_{ij}(t,x) y_i y_j \ge 0$, $x, y \in \mathbb{R}^d$, $0 \le t \le T$.
(P2) There exists a number $C > 0$ so that
$$\sum_{k \ge 1} \left( \sup_{t,x} |q_k \nu_k(t,x)|^2 + \int_0^T q_k^2 \|g_k\|_{L_2(\mathbb{R}^d)}^2(t) \, dt \right) \le C.$$

For the matrix $A$ and each $t, x$, we have $A_{ij}(t,x) = \tilde{\sigma}_{ik}(t,x) \tilde{\sigma}_{jk}(t,x)$, where the functions $\tilde{\sigma}_{ik}$ are bounded. This representation might not be unique; see, for example, [9], Theorem III.2.2 or [36], Lemma 5.2.1. Given any such representation of $A$, consider the following backward Itô equation:

$$X_{t,x,i}(s) = x_i + \int_s^t B_i(\tau, X_{t,x}(\tau)) \, d\tau + \sum_{k \ge 1} \int_s^t q_k \sigma_{ik}(\tau, X_{t,x}(\tau)) \overleftarrow{dw_k}(\tau)$$
(5.2)
$$+ \int_s^t \tilde{\sigma}_{ik}(\tau, X_{t,x}(\tau)) \overleftarrow{d\widetilde{w}_k}(\tau); \qquad s \in (0,t), t \in (0,T], t\text{-fixed},$$

where $B_i = b_i - \sum_{k\geq 1} q_k^2 \sigma_{ik}\nu_k$ and $\tilde{w}_k, k \geq 1$, are independent standard Wiener processes on $\mathbb{F}$ that are independent of $w_k, k \geq 1$. This equation might not have a strong solution, but does have weak, or martingale, solutions due to assumptions (B1)–(B3) and properties (P1) and (P2) of the sequence $Q$; this weak solution is unique in the sense of probability law ([36], Theorem 7.2.1).

Let $\mathcal{Q}$ be the operator

$$\mathcal{Q}:\{u_\alpha, \alpha \in \mathcal{J}\} \mapsto \{q^\alpha u_\alpha, \alpha \in \mathcal{J}\}.$$

THEOREM 5.1. *Under assumptions* (B0)–(B5), (5.1) *has a unique Wiener Chaos solution* $u = u(t,x)$. *If* $Q$ *is a sequence with properties* (P1) *and* (P2), *then* $u(t,\cdot) \in L_{2,Q}(\mathbb{W}; L_2(\mathbb{R}^d))$, $0 \leq t \leq T$, *and the following representation holds:*

$$
\begin{aligned}
u(t,x) = \mathcal{Q}^{-1}\mathbb{E}\bigg(&\int_0^t f(s, X_{t,x}(s))\gamma(t,s,x)\, ds \\
&+ \sum_{k\geq 1}\int_0^t q_k g_k(s, X_{t,x}(s))\gamma(t,s,x)\,\overleftarrow{dw_k}(s) \\
&+ u_0(X_{t,x}(0))\gamma(t,0,x)\big|\mathcal{F}_t^W\bigg),
\end{aligned}
$$
(5.3)

$t \leq T, x \in \mathbb{R}^d$, *where* $X_{t,x}(s)$ *is a weak solution of* (5.2) *and*

$$
\begin{aligned}
\gamma(t,s,x) = \exp\bigg(&\int_s^t c(\tau, X_{t,x}(\tau))\, d\tau + \sum_{k\geq 1}\int_s^t q_k \nu_k(\tau, X_{t,x}(\tau))\,\overleftarrow{dw_k}(\tau) \\
&- \tfrac{1}{2}\int_s^t \sum_{k\geq 1} q_k^2 |\nu_k(\tau, X_{t,x}(\tau))|^2\, d\tau\bigg).
\end{aligned}
$$
(5.4)

PROOF. It is enough to establish (5.3) when $t = T$. Consider the equation

$$
\begin{aligned}
dU = &(a_{ij}D_iD_jU + b_iD_iU + cU + f)\, dt \\
&+ \sum_{k\geq 1}(\sigma_{ik}D_iU + \nu_kU + g_k)q_k\, dw_k
\end{aligned}
$$
(5.5)

with initial condition $U(0,x) = u_0(x)$. In Theorem 4.1, take $Q = 1$ and consider the normal triple $(H_2^1(\mathbb{R}^d), L_2(\mathbb{R}^d), H_2^{-1}(\mathbb{R}^d))$. Then (5.5) has a unique Wiener Chaos solution $U$ and

$$U(t,\cdot) \in L_2(\mathbb{W}; L_2(\mathbb{R}^d)), \qquad 0 \leq t \leq T.$$

By construction, process $u = \mathcal{Q}^{-1}U$ is the corresponding Wiener Chaos solution of (5.1). To establish representation (5.3), consider the function $U_h = \mathbb{E}(U\mathcal{E}(h)), h \in \mathcal{H}$. According to Theorem 3.7, the function $U_h$ is the unique solution of the equation

$$
\begin{aligned}
dU_h &= (a_{ij}D_iD_jU_h + b_iD_iU_h + cU_h + f)\,dt \\
&\quad + \sum_{k\geq 1}(\sigma_{ik}D_iU_h + \nu_kU_h + g_k)q_kh_k\,dt
\end{aligned}
\tag{5.6}
$$

with initial condition $U_h|_{t=0} = u_0$. We also define

$$
\begin{aligned}
Y(T,x) &= \int_0^T f(s,X_{T,x}(s))\gamma(T,s,x)\,ds \\
&\quad + \sum_{k\geq 1}\int_0^T g_k(s,X_{T,x}(s))\gamma(T,s)q_k\,\overleftarrow{dw_k}(s) \\
&\quad + u_0(X_{T,x}(0))\gamma(T,0,x).
\end{aligned}
\tag{5.7}
$$

By direct computation,

$$\mathbb{E}(\mathbb{E}(\mathcal{E}(h)Y(T,x)|\mathcal{F}_T^W)) = \mathbb{E}(\mathcal{E}(h)Y(T,x)) = \mathbb{E}'Y(T,x),$$

where $\mathbb{E}'$ is the expectation with respect to the measure $d\mathbb{P}'_T = \mathcal{E}(h)\,d\mathbb{P}_T$ and $\mathbb{P}_T$ is the restriction of $\mathbb{P}$ to $\mathcal{F}_T^W$.

Under assumptions (B0)–(B5), the solution $U_h$ of (5.6) is continuous in $(t,x)$ and has a probabilistic representation via the Feynmann–Kac formula; see [21], Section 1.6. Using the Girsanov theorem ([14], Theorem 3.5.1), this representation can be written as $U_h(T,x) = \mathbb{E}'Y(T,x)$ or

$$\mathbb{E}(\mathcal{E}(h)U(T,x)) = \mathbb{E}(\mathcal{E}(h)\mathbb{E}Y(t,x)|\mathcal{F}_T^W).$$

By Remark 3.6, the last equality implies $U(T,x) = \mathbb{E}(Y(T,x)|\mathcal{F}_T^W)$ as elements of $L_2(\mathbb{W})$.

Theorem 5.1 is proved. $\square$

EXAMPLE 5.2 (Krylov–Veretennikov formula).   Consider the equation

$$du = (a_{ij}D_iD_ju + b_iD_iu)\,dt + \sum_{k=1}^d \sigma_{ik}D_iu\,dw_k, \qquad u(0,x) = u_0(x). \tag{5.8}$$

Assume (B0)–(B5) and suppose that $a_{ij}(t,x) = \frac{1}{2}\sigma_{ik}(t,x)\sigma_{jk}(t,x)$. By Corollary 4.2, (5.8) has a square integrable solution such that

$$\mathbb{E}\|u\|^2_{L_2(\mathbb{R}^d)}(t) \leq C^*\|u_0\|^2_{L_2(\mathbb{R}^d)}.$$

By Theorem 4.1 with $Q = 1$,

$$u(t,x) = \sum_{n=1}^\infty \sum_{|\alpha|=n} u_\alpha(t,x)\xi_\alpha$$

$$= u_0(x)$$

(5.9)
$$+ \sum_{n=1}^{\infty} \sum_{k_1,\ldots,k_n=1}^{d} \int_0^t \int_0^{s_n} \cdots \int_0^{s_2} P_{t,s_n} \sigma_{jk_n} D_j$$
$$\cdots P_{s_2,s_1} \sigma_{ik_1} D_i P_{s_1,0} u_0(x)\, dw_{k_1}(s_1) \cdots dw_{k_n}(s_n),$$

where $P_{t,s}$ is the semigroup generated by the operator $\mathcal{A} = a_{ij} D_i D_j u + b_i D_i u$. On the other hand, in this case, Theorem 5.1 yields

$$u(t,x) = \mathbb{E}(u_0(X_{t,x}(0))|\mathcal{F}_t^W),$$

where $W = (w_1, \ldots, w_d)$ and

(5.10)
$$X_{t,x,i}(s) = x_i + \int_s^t b_i(\tau, X_{t,x}(\tau))\, d\tau + \sum_{k=1}^{d} \int_s^t \sigma_{ik}(\tau, X_{t,x}(\tau)) \overleftarrow{dw_k}(\tau),$$

$$s \in (0,t), t \in (0,T], t\text{-fixed}.$$

Thus, we have arrived at the Krylov–Veretennikov formula (cf. [20], Theorem 4):

$$\mathbb{E}(u_0(X_{t,x}(0))|\mathcal{F}_t^W)$$

(5.11) $$= u_0(x) + \sum_{n=1}^{\infty} \sum_{k_1,\ldots,k_n=1}^{d} \int_0^t \int_0^{s_n} \cdots \int_0^{s_2} P_{t,s_n} \sigma_{jk_n} D_j$$
$$\cdots P_{s_2,s_1} \sigma_{ik_1} D_i u_{(0)}(x)\, dw_{k_1}(s_1) \cdots dw_{k_n}(s_n).$$

**6. Passive scalar in a Gaussian field.** The following viscous transport equation is used to describe time evolution of a passive scalar $\theta$ in a given velocity field $\mathbf{v}$:

(6.1) $\quad \dot{\theta}(t,x) = \nu \Delta \theta(t,x) - \mathbf{v}(t,x) \cdot \nabla \theta(t,x) + f(t,x); \qquad x \in \mathbb{R}^d, d > 1.$

In the Kraichnan model of turbulent transport [10, 18], using the results from [1] and [23], (6.1) can be written as an Itô stochastic evolution equation:

(6.2) $\quad d\theta(t,x) = (\nu \Delta \theta(t,x) + f(t,x))\, dt - \sigma_k(x) \cdot \nabla \theta(t,x)\, dw_k(t),$

where

(S1) $\{w_k(t), k \geq 1, t \geq 0\}$ is a collection of independent Wiener processes.
(S2) $\{\sigma_k, k \geq 1\}$ is an orthonormal basis in the space $H_C$, the reproducing kernel Hilbert space corresponding to the spatial covariance function $C$ of $\mathbf{v}$, so that $\sigma_k^i(x)\sigma_k^j(y) = \mathbb{E}v^i(x)v^j(y) = C^{ij}(x-y)$ and $C^{ij}(0) = \delta_{ij}$. The space $H_C$ is all or part of the Sobolev space $H^{(d+\gamma)/2}(\mathbb{R}^d; \mathbb{R}^d)$, $\gamma \in (0,2)$.

(S3) The initial condition $\theta_0$ is nonrandom and belongs to $L_2(\mathbb{R}^d)$.

Using the Wiener Chaos approach, it is possible to consider velocity fields **v** that are even more turbulent, for example,

$$v^i(t,x) = \sum_{k \geq 0} \sigma_k^i(x) \dot{w}_k(t), \tag{6.3}$$

where $\{\sigma_k, k \geq 1\}$ is an orthonormal basis in $L_2(\mathbb{R}^d; \mathbb{R}^d)$. With **v** as in (6.3), the passive scalar equation (6.2) becomes

$$\dot{\theta}(t,x) = \nu \Delta \theta(t,x) + f(t,x) - \nabla \theta(t,x) \cdot \dot{W}(t,x), \tag{6.4}$$

where $\dot{W} = \dot{W}(t,x)$ is a $d$-dimensional space–time white noise.

In [6] and [34], a similar equation is studied using the white noise approach. For more related results and references, see [12], Section 4.3. We will consider an even more general equation.

THEOREM 6.1. *Suppose that $\nu > 0$ is a real number, the functions $\sigma_k^i(x)$ are bounded and measurable, $\theta_0 \in L_2(\mathbb{R}^d)$ and $f \in L_2((0,T); H_2^{-1}(\mathbb{R}^d))$.*

*Fix $\varepsilon > 0$ and let $Q = \{q_k, k \geq 1\}$ be a sequence so that, for all $x, y \in \mathbb{R}^d$,*

$$2\nu |y|^2 - \sum_{k \geq 1} q_k^2 \sigma_k^i(x) \sigma_k^j(x) y_i y_j \geq \varepsilon |y|^2.$$

*Then there exists a unique Wiener Chaos solution of* (6.2). *This solution satisfies*

$$\theta \in L_{2,Q}(\mathbb{W}; L_2((0,T); H_2^1(\mathbb{R}^d))) \cap L_{2,Q}(\mathbb{W}; \mathbf{C}((0,T); L_2(\mathbb{R}^d))).$$

PROOF. The result follows from Theorem 4.1, with $\mathcal{A} = \nu \Delta$, $\mathcal{M}_k = \sigma_k^i D_i$, in the normal triple $(H_2^1(\mathbb{R}), L_2(\mathbb{R}), H_2^{-1}(\mathbb{R}))$. □

REMARK 6.2. If $\max_i \sup_x |\sigma_k^i(x)| \leq C_k$, $k \geq 1$, then a possible choice of $Q$ is $q_k = (\delta \nu)^{1/2}/(d 2^k C_k)$, $0 < \delta < 2$.

If $\sigma_k^i(x) \sigma_k^j(x) \leq C_\sigma < +\infty$, $i, j = 1, \ldots, d$, $x \in \mathbb{R}^d$, then a possible choice of $Q$ is $q_k = \varepsilon (2\nu/(C_\sigma d))^{1/2}$, $0 < \varepsilon < 1$.

When $\nu = 0$, (6.2) describes nonviscous transport [8] and can still be studied if interpreted in the Stratonovich sense:

$$du(t,x) = f(t,x) \, dt - \sigma_k(x) \cdot \nabla \theta(t,x) \circ dw_k(t). \tag{6.5}$$

In what follows, we assume that $f = 0$, each $\sigma_k$ is divergence free and $\sigma_k^i(x) \sigma_k^j(x) = \delta_{ij}$. Then (6.5) has an equivalent Itô form

$$d\theta(t,x) = \tfrac{1}{2} \Delta \theta(t,x) \, dt - \sigma_k^i(x) D_i \theta(t,x) \, dw_k(t). \tag{6.6}$$

The following result (cf. [25]) summarizes the main facts about the Wiener Chaos solution of (6.6)

THEOREM 6.3. *In addition to* (S1)–(S3), *assume that each $\sigma_k$ is divergence free. Then there exists a unique Wiener Chaos solution $\theta = \theta(t,x)$ of* (6.6). *This solution has the following properties:*

  (i) *For each $t \in [0,T]$, $\theta(t,\cdot) \in L_2(\mathbb{R}^d)$.*
  (ii) *For every $\varphi \in \mathbf{C}_0^\infty(\mathbb{R}^d)$ and all $t \in [0,T]$, the equality*

$$(6.7) \quad (\theta,\varphi)(t) = (\theta_0,\varphi) + \tfrac{1}{2}\int_0^t (\theta,\Delta\varphi)(s)\,ds + \int_0^t (\theta,\sigma_k^i D_i\varphi)\,dw_k(s)$$

*holds in $L_2(\mathbb{W})$, where $(\cdot,\cdot)$ is the inner product in $L_2(\mathbb{R}^d)$.*

  (iii) *If $\theta_0 \in W_p^2(\mathbb{R}^d)$, $p \geq d+1$, and $X$ is a weak solution of*

$$(6.8) \quad X_{t,x} = x + \int_0^t \sigma_k(X_{s,x})\,dw_k(s),$$

*then*

$$(6.9) \quad \theta(t,x) = \mathbb{E}(\theta_0(X_{t,x})|\mathcal{F}_t^W), \qquad 0 \leq t \leq T, x \in \mathbb{R}^d.$$

  (iv) *For $1 \leq p < \infty$ and $r \in \mathbb{R}$, denote by $L_{p,(r)}(\mathbb{R}^d)$ the space of measurable functions with finite norm*

$$\|f\|_{L_{p,(r)}(\mathbb{R}^d)}^p := \int_{\mathbb{R}^d} |f(x)|^p (1+|x|^2)^{pr/2}\,dx.$$

*Then there exists a number $K$ depending only on $p$ and $r$ so that*

$$(6.10) \quad \mathbb{E}\|\theta\|_{L_{p,(r)}(\mathbb{R}^d)}^p(t) \leq e^{Kt}\|\theta_0\|_{L_{p,(r)}(\mathbb{R}^d)}^p, \qquad t > 0.$$

*In particular, if $r = 0$, then $K = 0$.*

PROOF. (i) This follows from Theorem 4.1, because (6.6) is weakly parabolic and one can take $Q = 1$.

(ii) Since for each $k$ we have $D_i\sigma_k^i = 0$ in the sense of distributions, the same arguments as in the proof of the second part of Theorem 3.4 result in (6.7).

(iii) Equality (6.9) follows from Theorem 5.1 after observing that the time-homogeneity of $\sigma_k$ allows us to rewrite the corresponding backward equation as (6.8).

(iv) To prove (6.10), denote by $S_t: \theta_0(\cdot) \mapsto \theta(t,\cdot)$, $t > 0$, the solution operator for (6.6). Direct calculations show that, for every $r \in \mathbb{R}$, assumptions (A1)–(A2) of Theorem 4.1 are satisfied with $V = H_{2,(r)}^1(\mathbb{R}^d) := \{f : f, D_i f \in L_{p,(r)}(\mathbb{R}^d)\}$, $H = L_{2,(r)}(\mathbb{R}^d)$, $V' = H_{2,(r)}^{-1}(\mathbb{R}^d)$ and $Q = 1$. Then $S_t$ is a bounded linear operator from $L_{2,(r)}(\mathbb{R}^d)$ to $L_{2,(r)}(\Omega \times \mathbb{R}^d)$ and

$$\|S_t\theta_0\|_{L_{2,(r)}(\Omega \times \mathbb{R}^d)} \leq e^{K_2 t}\|\theta_0\|_{L_{2,(r)}(\mathbb{R}^d)}$$

for every $r \in \mathbb{R}$, where the number $K_2$ depends only on $r$, and

$$\|f\|^2_{L_{2,(r)}(\Omega \times \mathbb{R}^d)} := \mathbb{E}\|f\|^2_{L_{2,(r)}(\mathbb{R}^d)}.$$

The calculations also show that $K_2 = 0$ if $r = 0$.

By (6.9), $S_t$ is a bounded linear operator from $L_\infty(\mathbb{R}^d)$ to $L_\infty(\Omega \times \mathbb{R}^d)$ and

$$\|S_t \theta_0\|_{L_\infty(\Omega \times \mathbb{R}^d)} \leq \|\theta_0\|_{L_\infty(\mathbb{R}^d)}.$$

Interpreting the last inequality in the weighted spaces $L_{\infty,(r)}$ with weight $r = 0$, we interpolate between $L_\infty = L_{\infty,(0)}$ and $L_{2,(r)}$ and derive (6.10) from the Stein interpolation theorem in weighted spaces (see, e.g., [2], Theorem 4.3.6). Theorem 6.3 is proved. $\square$

**Acknowledgment.** We are indebted to the referee for a number of insightful suggestions.

DEPARTMENT OF MATHEMATICS
UNIVERSITY OF SOUTHERN CALIFORNIA
3620 S. VERMONT AVE., KAP 108
LOS ANGELES, CALIFORNIA 90089-2532
USA
E-MAIL: lototsky@math.usc.edu
   rozovski@math.usc.edu
URL: www.usc.edu/schools/college/mathematics/people/faculty/lototsky.html
   www.usc.edu/dept/LAS/CAMS/usr/facmemb/boris/main.htm